\documentstyle{amsppt}
\NoRunningHeads

\def\r2{\Bbb R^2}
\def\z2{\Bbb Z^2}
\def\z+d{\Bbb Z_+^d}
\def\rd{\Bbb R^d}

\def\ve{\varepsilon}

\def\ol{\overline l}
\def\og{\overline \gamma}

\def\dim{\hbox{dim}}

\def\rr{\Bbb R}

\topmatter
\author{Sergei Konyagin and Izabella {\L}aba}
\endauthor
\title
Separated sets and the Falconer conjecture for polygonal norms
\endtitle
\address Department oF Mechanics and Mathematics, Moscow State University,
Moscow, 119992, Russia,
e-mail: konyagin\@ok.ru
\endaddress

\address Department of Mathematics, University of British Columbia,
Vancouver, B.C. V6T 1Z2, Canada, e-mail: ilaba\@math.ubc.ca
\endaddress
\abstract
The Falconer conjecture \cite{F86} asserts that if $E$ is a planar set with Hausdorff
dimension strictly greater than 1, then its Euclidean distance set $\Delta(E)$
has positive one-dimensional Lebesgue measure.  We discuss the analogous
question with the Euclidean distance replaced by non-Euclidean norms $\|\cdot\|_X$
in which the unit ball is a polygon with $2K$ sides.  We prove that
for any such norm, and for any $\alpha>K/(K-1)$, there is a set of
Hausdorff dimension $\alpha$ whose
distance set has Lebesgue measure 0.

Mathematics Subject Classification: 28A78.

\endabstract
\endtopmatter
\document

\centerline{\S0. INTRODUCTION}
\bigskip

A conjecture of Falconer \cite{F86} asserts that if a set $E\subset\rr^2$ has Hausdorff
dimension strictly greater than 1, then its Euclidean distance set
$$\Delta(E)=\Delta_{l_2^2}(E)=\left\{\|x-x'\|_{l_2^2}:x,x'\in E\right\}$$
has positive one-dimensional Lebesgue measure.  The current best result
in this direction is due to Wolff \cite{W99}, who proved that the conclusion
is true if $E$ has Hausdorff dimension greater than 4/3.  Erdogan \cite{Er03}
extended this result to higher dimensions, proving that the same conclusion
holds for subsets of $\rr^d$ with Hausdorff dimension greater than
$d(d+2)/2(d+1)$.  This improves
on the earlier results of Falconer \cite{F86}, Mattila \cite{M87},
and Bourgain \cite{B94}.

A similar question can be posed for more general two-dimensional normed
spaces. More precisely, if $X$ is such a space and $E\subset X$, then
we define the $X$-distance set of $A$ as
$$\Delta_{X}(E)=\left\{\|x-x'\|_X:x,x'\in E\right\}$$
and ask how the size of $\Delta_{X}(E)$ depends on the dimension
of $E$ as well as on the properties of the norm $\|\cdot\|_X$.
Simple examples show that Falconer's conjecture as stated above,
but with $\Delta(E)$ replaced by $\Delta_X(E)$, cannot hold for
all normed spaces $X$.  For instance, let
$$\|x\|_{l_\infty^2}=\max(|x_1|,|x_2|)$$
and let $E=F\times F$, where $F$ is a subset of $[0,1]$
with Hausdorff dimension 1 such that $F-F:=\{x-x':x,x'\in F\}$ has
measure 0.  (It is an easy exercise to modify the Cantor set construction
to produce such a set.)  Then $E$ has Hausdorff dimension 2, but
its $l_\infty^2$-distance set $F-F$ has measure 0.

Here and below, we use $\dim(E)$ to denote the Hausdorff dimension
of $E$, $|F|_d$ to denote the $d$-dimensional Lebesgue measure
of $F$, and $|A|$ to denote the cardinality of a finite set $A$.

\proclaim{Definition 0.1.}  Let $0<\alpha<2$.  We will say that
the $\alpha$-Falconer conjecture holds in $X$ if for any set
$E\subset X$ with $\dim(E)>\alpha$ we have $|\Delta_X(E)|_1>0$.
\endproclaim

Iosevich and the second author \cite{I{\L}04} proved that the
$3/2$-Falconer conjecture holds if the unit ball in $X$,
$$BX=\{x\in\r2:\,\|x\|_X\le1\},$$
is strictly convex and its boundary $\partial BX$ has everywhere
nonvanishing curvature, in the sense that the diameter of the
chord
$$\{x\in BX:x\cdot v\geq\max_{y\in BX}(y\cdot v)-\epsilon\},$$
where $v$ is a unit vector and $\epsilon>0$, is bounded by
$C\sqrt{\epsilon}$ uniformly for all $v$ and $\epsilon$.
We do not know of any counterexamples to the 1-Falconer conjecture
in normed spaces with $BX$ strictly convex.

On the other hand, if $BX$ is a polygon, then the above example
shows that the $\alpha$-Falconer conjecture may fail for all
$\alpha <2$.  The purpose of this paper is to examine this
situation in more detail.

\proclaim{Theorem 1} Let $BX$ be a symmetric convex polygon with $2K$
sides.  Then there is a set $E\subset
[0,1]^2$ with Hausdorff dimension $\ge K/(K-1)$ such that
$|\Delta_X(E)|_1=0$.
\endproclaim

If we assume that there is a coordinate system in which
the slopes of all sides of $K$ are algebraic, then a stronger result
is known \cite{K{\L}04}.

\proclaim{Corollary 2} \cite{K{\L}04} If $BX$ is a polygon with
finitely many sides, and if there is a coordinate system in which
all sides of $BX$ have algebraic slopes,
then there is a compact $E\subset X$ such that the Hausdorff dimension of $E$
is $2$ and the Lebesgue measure of $\Delta_X(E)$ is $0$.
\endproclaim

In particular, Corollary 2 can be applied to all polygons $BX$
with $4$ or $6$ sides. We do not know if the same assertion is
true for all polygonal norms. However, using recent results on
Diophantine approximations, one can prove it for almost all polygons
$BX$. Fixing a coordinate system, we can define,
for any non-degenerate segment $I\subset X$, its slope $Sl(I)$: if the line
containing $I$ is given by an equation $u_1x_1+u_2x_2+u_0=0$, then
we set $Sl(I)=-u_1/u_2$. We write $Sl(I)=\infty$ if $u_2=0$.

\proclaim{Theorem 3} For any integer $K\ge2$ and for almost all
$\gamma_1,\dots,\gamma_K$ the following is true. If $BX$ is
a symmetric convex polygon with $2K$ sides, and the slopes of
non-parallel sides are equal to $\gamma_1,\dots,\gamma_K$, then
there is a compact $E\subset X$ such that the Hausdorff dimension of $E$
is $2$ and the Lebesgue measure of $\Delta_X(E)$ is $0$.
\endproclaim

Actually, we will prove the stronger result: if the slopes of $3$ non-parallel
sides of $BX$ are fixed, then for almost all choices of slopes of
other $K-3$ non-parallel sides the required compact $A$ exists
(recall that for $K\le3$ Theorem 3 follows from Corollary 2).

\newpage

\centerline{\S1. PROOF OF THEOREM 1}
\bigskip

We may assume that $K\geq 4$, since otherwise Corollary 2 applies.
We use $B(x,r)$ to denote the closed Euclidean ball with center at $x$ and
with radius $r$.  We also denote $A-A=\{a-a':\ a,a'\in A\}$ and
$A\cdot v=\{a\cdot v:\ a\in A\}$.

Let $b_1,\dots,b_K$ be vectors such that
$$BX=\bigcap_{k=1}^K\{x:\ |x\cdot b_k|\leq 1\}.$$
Then for any $x\in X$,
$$\|x\|_X=\max_{1\leq k\leq K}|x\cdot b_k|.\tag1.3$$
Let also $a_1,\dots,a_K$ be unit vectors parallel to the $K$ sides
of $BX$, so that
$$a_j\cdot b_j=0,\ j=1,\dots,K.\tag1.4$$

\proclaim{Lemma 1.1}
Assume that $K\geq 4$.  Then there are arbitrarily large integers $n$
for which we may choose sets $A=A(n)\subset B(0,1/2)$ such that $|A|=n$ and
$$|(A-A)\cdot b_k|\ll n^{1-1/K},\ k=1,2,\dots,K,\tag1.1$$
(in particular, $|\Delta_X(A)|\ll n^{1-1/K}$), and
$$\|x-x'\|_X\gg n^{-1/2},\ x,x'\in A, \ x\neq x',\tag1.2$$
with the implicit constants independent of $n$.
\endproclaim

\demo{Proof}
Fix a large integer $N$, and let $u_1,\dots,u_K$
be numbers in $[1,2]$, to be determined later.  Define
$$S=\Big\{\sum_{k=1}^K \frac{j_k}{N}u_ka_k,\ j_k\in\{1,\dots,N\}\Big\}.$$
We claim that the set
$$U=\{(u_1,\dots,u_K)\in \rr^{K}:\ |S|<N^K\}$$
has $K$-dimensional measure 0.  Indeed, if $|S|<N^K$, then we must have
$$
\sum_{k=1}^K \frac{j_k}{N}u_ka_k=0$$
for some $j_1,\dots,j_K\in\{1-N,\dots,N-1\}$, not all zero.  Fix such
$j_1,\dots,j_K$.  Then
the $2\times K$ matrix with columns $\frac{j_k}{N}u_ka_k$, $k=1,\dots,K$,
has rank at least 1, hence its nullspace has dimension at most $K-1$.
It follows that $U$ is a union of a finite number of hyperplanes of
dimension at most $K-1$, therefore has $K$-dimensional measure $0$
as claimed.

We will assume henceforth that $(u_1,\dots,u_K)\notin U$.
Then $|S|=N^K$ and $S\subset B(0,2K)$.  Our goal is to
obtain (1.1), (1.2) for $n=N^K$ and $A=(4K)^{-1}S$.

We first prove that (1.1) holds, i.e.
$$|(S-S)\cdot b_k| \ll N^{K-1}\ll n^{1-1/K},\ k=1,2,\dots,K.\tag1.5$$
Indeed, let $x\in S-S$, then $x=\sum_{k=1}^K \frac{j_k}{N}u_ka_k$
for some $j_1,\dots,j_K\in\{1-N,\dots,N-1\}.$
Fix $k_0\in\{1,\dots,k\}$, then
$$x\cdot b_{k_0}=\sum_{k=1}^K \frac{j_k}{N}u_ka_k\cdot b_{k_0}
=\sum_{k\neq k_0} \frac{j_k}{N}u_ka_k\cdot b_{k_0},$$
where we also used (1.4).  The last sum can take at most
$(2N)^{K-1}$ possible values, which proves (1.5).

It remains to verify that there is a choice of $u_1,\dots,u_K$
for which (1.2) also holds.  We will do so by proving that
if $t$ is a sufficiently small constant, depending only on $K$ and
on the angles between the non-parallel sides of $BX$, then the set
$$\{(u_1,\dots,u_K)\in [1,2]^K:\ \|x\|_X\leq tN^{-K/2}
\hbox{ for some } x\in S-S\}\tag1.6$$
has $K$-dimensional Lebesgue measure strictly less than 1.

Let $x\in S-S$, then $x=\sum_{k=1}^K \frac{j_k}{N}u_ka_k$
for some $j_k\in\{1-N,\dots,N-1\}.$  Suppose that $x\neq 0$ and
$$\|x\|_X\leq tN^{-K/2}.\tag1.7$$
Assume that $|j_{k_1}|\geq |j_{k_2}|\geq\dots\geq |j_{k_K}|$,
and that $|j_{k_1}|\in [2^s,2^{s+1})$ for some integer $s$ such that
$1\leq 2^s\leq N$.  If we had $|j_{k_2}|<2^{s-2}/K$, then we would also
have
$$
\|x\|_X\geq
\|\frac{j_{k_1}}{N}u_{k_1}a_{k_1}\|_X
-\sum_{k\neq k_1} \|\frac{j_k}{N}u_ka_k\|_X
\geq\frac{2^s}{N}-K\cdot \frac{2\cdot 2^{s-2}}{KN}=
\frac{2^{s-1}}{N}\geq \frac{1}{2N}.$$
But if $K\geq 4$, then (1.7) implies that $\|x\|_X\leq tN^{-2}$,
which contradicts the last inequality if $t\leq 1$ and $N>2$.
It follows that
$$|j_{k_1}|\geq 2^s,\ |j_{k_2}|\geq 2^{s-2}/K.\tag1.8$$

Fix $j_{k_1},j_{k_2}$ as in (1.8).
Fix also $y=\sum_{k\neq k_1,k_2} \frac{j_k}{N}u_ka_k$, and consider
the set of $(u_{k_1},u_{k_2})\in\rr^2$ such that (1.7) holds,
i.e.
$$\|\frac{j_{k_1}}{N}u_{k_1}a_{k_1}
+\frac{j_{k_2}}{N}u_{k_2}a_{k_2}+y\|_X\leq tN^{-K/2}.$$
By (1.8), this set has 2-dimensional measure
$$\leq
c_1(tN^{-K/2})^2\cdot\frac{N}{2^s}\cdot\frac{NK}{2^{s-2}}
=4c_1K\cdot t^2N^{2-K}/2^{2s}.$$
Here and through the rest of the proof of the lemma, $c_1,c_2,c_3$
denote constants which may depend on $K$ and on the angles between
the non-parallel sides of $BX$, but are independent of $t$ and $N$.

Integrating over $u_k$, $k\neq k_1,k_2$, we see that the set
$$\Big\{(u_1,\dots,u_K)\in [1,2]^K:\
\|\sum_{k=1}^K \frac{j_k}{N}u_ka_k\|_X\leq tN^{-K/2}\Big\},$$
with fixed $j_1,\dots,j_K$ such that
$$2^s\leq\max_{k=1,\dots,K}|j_k|<2^{s+1},\tag1.9$$
has $K$-dimensional measure $\leq 4c_1K\cdot t^2 N^{2-K}/2^{2s}$.

The number of $K$-tuples $j_1,\dots,j_K$ satisfying (1.9)
is $\leq (2^{s+2})^K$, hence summing over all such $K$-tuples we
get a set of measure
$$\leq c_2 t^2N^{2-K}2^{(K-2)s}.$$
Now sum over all $s$ with $2^s\leq N$.  We find that the measure of
the set in (1.6) is
$$\leq c_2 \sum_{s:1\leq 2^s\leq N}t^2N^{2-K}2^{(K-2)s}
\leq c_3 t^2 N^{2-K}N^{K-2}=c_3t^2.$$
This is less than 1 if $t<\sqrt{c_3}$, as claimed.
\enddemo

\demo{Proof of Theorem 1}
We construct $E$ as follows.  Take a small positive number $c$
which will be specified later. Let $A_j=A(n_j)$ be as in Lemma 1.1,
where a nondecreasing sequence $\{n_j\}$ and a sequence $\{N_j\}$ are
such that
$$N_j=\prod_{\nu=1}^{j} n_\nu,\quad n_j\to\infty\,(j\to\infty),
\quad \log n_{j+1}/\log N_j\to0\,(j\to\infty).\tag 1.10$$
(We consider that the empty product for $j=0$ is equal to $1$.)
Also, fix $s=(K-1)/K>1/2$.  Let also $c$ be small enough
so that for any $j$ the discs $B(x,cn_j^{-s})$, $x\in A_j$, are mutually
disjoint and contained in $B(0,1)$; this is possible by (1.2).
Denote
$$\delta_j=cn_j^{-s},\quad \Delta_j=\prod_{\nu=1}^{j} \delta_j=c^j N_j^{-s}.$$
Let $E_1=\bigcup_{x\in A_1}B(x,\delta_1)$.
We then define $E_2,E_3,\dots$ by induction.
Namely, suppose that we have constructed $E_j$ which is
a union of $N_j$ disjoint closed discs $B_i$ of radius $\Delta_j$ each.
Then $E_{j+1}$ is obtained from $E_j$ by replacing each $B_i$ by
the image of $\bigcup_{x\in A_{j+1}}B(x,\delta_{j+1})$ under the unique affine mapping
which takes $B(0,1)$ to $B_i$ and preserves direction of vectors.
We then let $E= \bigcap_{j=1}^\infty E_j$.

We will first prove that $E$ has Hausdorff dimension at least $1/s$.  The
calculation follows closely that in \cite{F85}, pp. 16--18.

Let ${\Cal B}_j$ be the family of all discs of radius $\Delta_j$
used in the construction of $E_j$,
and let ${\Cal B}=\bigcup_{ j=0}^\infty {\Cal B}_j$,
where we set ${\Cal B}_0=\{B(0,1)\}$.
We then define
$$\mu(F)=\inf\Big\{\sum_{i=1}^\infty N_{j(i)}^{-1}:\ F\subset
\bigcup_{i=1}^\infty B(x_i,r_i),
B(x_i,r_i)\in{\Cal B_{j(i)}}\Big\},\tag1.11$$
for all $F\subset E$.  Clearly, $\mu$ is an outer measure on
subsets of $E$. Observe that if $B=B(x,\Delta_j)\in{\Cal B}_j$, then
$$N_j^{-1}=n_{j+1}\cdot N_{j+1}^{-1}
=\sum_{B'\in{\Cal B}_{j+1}:B'\subset B}(N_{j+1})^{-1},\tag1.12$$
hence the sum in (1.11) does not change if we replace a disc $B
\in{\Cal B}_j$ by all its subdiscs from the next iteration ${\Cal B}_{j+1}$.
In particular, we may assume that all the discs in the covering of
$F$ in (1.11) have radius less than $\delta$ for any $\delta>0$.

We first claim that if $B_0=B_0(x_0,r_0)\in{\Cal B}_j$ then
$$\mu(E\cap B_0)=N_j^{-1}.\tag1.13$$
The inequality $\mu(E\cap B_0)\leq N_j^{-1}$ is obvious, by taking a
covering of $E\cap B_0$ by the single ball $B_0$.
Let now $E\cap B_0\subset\bigcup_i B_i$, where $B_i\in{\Cal B}$
has radius $r_i=\Delta_{j(i)}$. We need to prove that
$$\sum r_i^{1/s}\geq r_0^{1/s}.\tag1.14$$
Since $E$ is compact and $B_i$ are open relative
to $E$, we may assume that the covering is finite.  We may also assume
that all $B_i$ are disjoint, since otherwise we may simply remove any
discs contained in any other disc of the covering.  If the covering
consists of the single disc $B_0$, we are done.  Otherwise, let $B_I$ be
one of the covering discs with smallest $r_i$, say $B_I\in {\Cal B}_j$,
and let $\tilde B_I\in{\Cal B}_{j-1}$ be such that $B_I\subset \tilde B_I$.
Then $\tilde B_I\subset B_0$, hence all discs in ${\Cal B}_j$ contained
in $\tilde B_I$ are also contained in $B_0$.  By the minimality of
$r_I$, these discs belong to the covering $\{B_i\}$.  We then replace
all these discs by the single disc $\tilde B_I$; by (1.12), the
sum on the left side of (1.14) does not change.  Iterating this
procedure, we eventually arrive at a covering consisting only of
$B_0$, which proves (1.14).

Next, we prove that for any $s'>s$
$$\mu(E\cap B)\ll r^{1/s'}\tag1.15$$
for any disc $B=B(x,r)$, not necessarily in ${\Cal B}$,
where the constant in $\ll$ may depend on $s'$.
We may assume that $r\leq 1$, since otherwise we have from (1.13)
with $B_0=B(0,1)$
$$\mu(E\cap B)\leq \mu(E)=1\leq r^{1/s'},$$
which proves (1.15).  Let $j\geq 0$ be such that
$r\in(\Delta_{j+1},\Delta_j]$, and consider all discs in ${\Cal B}_j$
which intersect $E\cap B$.  They are closed, mutually disjoint discs
which intersect $B$ and have radius no less than $r$; hence there
are at most 6 such discs.  Applying (1.13) to each of these discs
and summing up, we have
$$\mu(E\cap B)\leq 6N_j^{-1}.$$
Moreover,
$$r>\Delta_{j+1}=N_j^{-s}n_{j+1}^{-s}c^{-j-1},$$
and we get (1.15) using (1.10).

Thus, if $s'>s$ and $\{B_i\}_{i=1}^\infty$ is a covering of $E$
by discs of radii $r_i$, then from (1.15) we have
$$\sum_{i=1}^\infty r_i^{1/s'}\gg
\sum_{i=1}^\infty \mu(E\cap B_i)\ge\mu(E).$$
Taking the infimum over all such coverings, we see that
$$H_{1/s'}(E)>0.$$
Since $s'>s$ is arbitrary, we conclude that
the Hausdorff dimension of $E$ is at least $K/(K-1)$.

It remains to prove that
$|\Delta_X(E)|_1=0$. From (1.1) we have
$$|(A-A)\cdot b_k|\leq Cn^{1-1/K},\ k=1,2,\dots,K,\tag1.16$$
with $C$ independent of $n$. We choose $c$ small enough so that
$$cC<1/2.\tag1.17$$

Let $D_j$ be the set of the centers of the discs in ${\Cal B}_j$.  We
claim that
$$|(D_j-D_j)\cdot b_k|\leq C^jN_j^s,\ k=1,2,\dots,K.\tag1.18$$
Indeed, for $j=1$ this is (1.16).  Assuming (1.18) for $j$, we now
prove it for $j+1$.  Let $x,x'\in D_{j+1}$.  Then $x\in B(y,\Delta_j)$,
$x'\in B(y',\Delta_j)$, $y,y'\in D_j$.  We write
$$(x-x')\cdot b_k=(y-y')\cdot b_k+((x-y)-(x'-y'))\cdot b_k.\tag1.19$$
The first term on the right is in $(D_j-D_j)\cdot b_k$, hence has
at most $C^jN_j^s$ possible values.  Also, by construction
$x-y,x'-y'$ are in $\Delta_jA_{j+1}$, hence the second term is in
$\Delta_j(A_{j+1}-A_{j+1})\cdot b_k$ and has at most $Cn_{j+1}^s$
possible values, by (1.16).  This gives at most $C^{j+1}N_{j+1}^s$ possible
values for (1.19), as required.

By (1.18), (1.3) and the triangle inequality, $\Delta_X(E_j)$ can be covered
by at most $KC^jN_j^s$ intervals of length $2c_0\Delta_j=2c_0c^jN_j^{-s}$,
where $c_0$ is the $X$-diameter of $B(0,1)$.  It follows that
$$|\Delta_X(E_j)|_1\leq 2Kc_0(cC)^j \leq 2Kc_0(1/2)^j,$$
by (1.17).  The last quantity goes to 0 as $j\to\infty$.
Since $\Delta_X(E)\subset\Delta_X(E_j)$, this proves our claim
that $|\Delta_X(E)|_1=0$.  The proof of the theorem is complete.

\enddemo

Remark. It is easy to check that the set constructed in the proof of Theorem 1
has the Hausdorff dimension exactly $K/(K-1)$.

\bigskip

\centerline{\S2. PROOF OF THEOREM 3}
\bigskip

The case $K\le3$ is covered by Corollary 2. We consider that $K>3$
and denote $d=K-3$. Denote
$$\ol=(l_1,\dots,l_d)\in\z+d,$$
$$\Cal L(L)=\{\ol:\,0\le l_k<L\,(k=1,\dots,d)\}.$$
For a real vector $\og=(\gamma_1,\dots,\gamma_d)$ we write
$\og\in(KM)$ if for any positive integer $L$ and for any $\ve>0$
$$\inf\left|\sum_{\ol\in \Cal L(L)} n_{\ol}
\gamma_1^{l_1}\dots\gamma_d^{l_d}\right|
\left(\max_{\ol\in \Cal L(L)}|n_{\ol}|\right)^{(1+\ve)L^d}>0,$$
where infimum is taken over all nonzero integral vectors
$\{n_{\ol}:\,\ol\in\Cal L\}$. The following theorem easily follows
from the results of Kleinbock and Margulis [KM98].

\proclaim{Theorem A} For almost all $\og\in\rd$ we have $\og\in (KM)$.
\endproclaim

The results of [KM98] have been refined in [BKM01], [Be02], [BBKM02].

Now we formulate the main result of this section.

\proclaim{Theorem 4} Let $\og\in (KM)$, $K=d+3$, and let $BX$ be
a symmetric convex polygon with $2K$ sides, and the slopes of
non-parallel sides are equal to $\gamma_1,\dots,\gamma_d,0,1$,
and $\infty$, then there is a compact $E\subset X$ such that the
Hausdorff dimension of $E$ is $2$ and the Lebesgue measure of
$\Delta_X(E)$ is $0$.
\endproclaim

Formally, Theorem 4 deals with polygons $BX$ of special kind, but it is easy
to see that for any polygon we can make slopes of three sides of it equal to
$0,1,\infty$ by a choice of a coordinate system. Indeed, if $I_1,I_2,I_3$
are 3 non-parallel sides of $BX$, then, taking the $x_1$-coordinate axis and
the $x_2$-coordinate axis of a new coordinate system parallel to
$I_1$ and $I_3$ respectively, we get $Sl(I_1)=0$,
$Sl(I_3)=\infty$; moreover, the slope of $I_2$ can be made equal to $1$ by
scaling and, if necessary, reflecting, the $x_2$-coordinate axis. Thus, combining
Theorem A and Theorem 4 we get Theorem 3 (and also its stronger version mentioned
in the end of \S0).

We use notation introduced in the beginning of \S1. To prove Theorem 4,
we need a lemma similar to Lemma 1.1.

\proclaim{Lemma 2.1}
Assume that $K, d, \og, BX$ satisfy the conditions of Theorem 4.
Then for any $\ve>0$ there are arbitrarily large integers $n$
for which we may choose sets $A=A(n)\subset B(0,1/2)$ such that $|A|=n$ and
$$|(A-A)\cdot b_k|\ll n^{(1/2)+\ve},\ k=1,2,\dots,K,\tag2.1$$
(in particular, $|\Delta_X(A)|\ll n^{(1/2)+\ve}$), and
$$\|x-x'\|_X\gg n^{-1/2-\ve},\ x,x'\in A, \ x\neq x',\tag2.2$$
where the implicit constants may depend on $\epsilon$ but are independent of $n$.
\endproclaim

\demo{Proof}
Fix a positive integer $L>1/\ve$. Next,
fix a large integer $N$. Define
$$S_0=\left\{\sum_{\ol\in \Cal L(L)}\frac{j_{\ol}}{N}\gamma_1^{l_1}\dots
\gamma_d^{l_d}:\,j_{\ol}\in\{1,\dots,N\}\right\}.\tag2.3$$
and $S=S_0\times S_0$, that is
$$S=\{(x_1,x_2):\,x_1,x_2\in S_0\}.$$
For any $x\in S_0$ we have
$$|x|\le\sum_{\ol\in \Cal L(L)}|\gamma_1|^{l_1}\dots|\gamma_d|^{l_d}
=\sum_{l=0}^{L-1}|\gamma_1|^l\dots\sum_{l=0}^{L-1}|\gamma_d|^l
\le\gamma^{dL},$$
where
$$\gamma=\max(|\gamma_1|,\dots,|\gamma|_d)+1.$$
Therefore, $S\subset B(0,2\gamma^{dL})$.  Our goal is to
check that $|S|=n$ and to obtain (2.1), (2.2) for $n=N^{2L^d}$ and
$A=(4\gamma^{dL})^{-1}S$.

We consider that $a_k$ ($k=1,\dots,d$) are parallel to the sides with slopes
$\gamma_1,\dots,\gamma_d$ respectively and $a_{d+1},a_{d+2},a_{d+3}$
are parallel to the sides with slopes $0,1,\infty$ respectively. Thus, we can
take $b_k=(-\gamma_k,1)$ for $k=1,\dots,d$, $b_{d+1}=(0,1)$, $b_{d+2}=(-1,1)$,
$b_{d+3}=(1,0)$.

We first prove (2.1) for $k=1,\dots,d$, i.e.
$$|(S-S)\cdot b_k|\ll n^{(1/2)+\ve}.\tag2.4$$
Indeed, for $x\in (S-S)\cdot b_{k_0}$, $k_0=1,2,\dots,d,$ we have
a representation
$$x\cdot b_{k_0}=-\gamma_k\sum_{\ol\in \Cal L(L)}\frac{j'_{\ol}}{N}
\gamma_1^{l_1}\dots\gamma_d^{l_d}
+\sum_{\ol\in \Cal L(L)}\frac{j''_{\ol}}{N}
\gamma_1^{l_1}\dots\gamma_d^{l_d},$$
where
$$j'_{\ol},j''_{\ol}\in\{1-N,\dots,N-1\}\quad(\ol\in\Cal L(L)).$$
Denote
$$\Cal L(L,k_0)=\{\ol:\,0\le l_k<L\,(k=1,\dots,d; k\neq k_0),\,
0\le l_{k_0}\le L\}.$$
Then we have
$$x\cdot b_{k_0}=\sum_{\ol\in \Cal L(L,k_0)}\frac{j_{\ol}}{N}
\gamma_1^{l_1}\dots\gamma_d^{l_d}$$
with
$$j_{\ol}\in\{2-2N,\dots,2N-2\}\quad(\ol\in\Cal L(L,k_0)).$$
Hence,
$$|(S-S)\cdot b_{k_0}|\ll (4N)^{L^d+L^{d-1}}.$$
By the choice of $L$ we have $L^d+L^{d-1}<(1+\ve)L^d$, and we get (2.4).
for $k=1,\dots,d$. Next, (2.4) holds for
$k=d+1,d+2,d+3$ because for those $k$ and for $x\in (S-S)\cdot b_k$
we have a representation
$$x\cdot b_{k_0}=\sum_{\ol\in \Cal L(L)}\frac{j_{\ol}}{N}
\gamma_1^{l_1}\dots\gamma_d^{l_d}$$
with
$$j_{\ol}\in\{2-2N,\dots,2N-2\}\quad(\ol\in\Cal L(L)).$$
Hence,
$$|(S-S)\cdot b_{k_0}|\le (4N)^{L^d},$$
and we again get (4.2) for sufficiently large $N$. So,
(2.1) is proved.

Now observe that the supposition $\og\in (KM)$ implies that
elements of $S_0$ with different representations (2.3) are distinct.
This gives $|S_0|=N^{L^d}$ and thus $|S|=|S_0|^2=n$ as required.
Moreover, since for any $x,x'\in S_0$ there is a representation
$$x-x'=\sum_{\ol\in \Cal L(L,k_0)}\frac{j_{\ol}}{N}
\gamma_1^{l_1}\dots\gamma_d^{l_d}$$
with
$$j_{\ol}\in\{1-N,\dots,N-1\}\quad(\ol\in\Cal L(L,k_0)).$$
we conclude from the supposition $\og\in (KM)$ that for $x\neq x'$
$$|x-x'|\gg (2N)^{-(1+0.1\ve)L^d-1}.\tag2.5$$
By the choice of $L$, we have $(1+0.1\ve)L^d+1\le(1+1.1\ve)L^d$,
and from (2.5) we get for sufficiently large $N$ and distinct $y,y'\in A$
$$\|y-y'\|_X\gg (4\gamma^{dL})^{-1}(2N)^{-(1+1.1\ve)L^d}
\gg N^{-(1+2\ve)L^d}=n^{-1/2-\ve}.$$
This completes the proof of Lemma 2.1.
\enddemo

\demo{Proof of Theorem 4}
We construct $E$ as follows. Let $A_j=A(n_j)$ be as in Lemma 2.1
with $\ve=\ve_j$, where a nondecreasing sequence $\{n_j\}$, a sequence
$\{N_j\}$, and a sequence $\{\ve_j\}$ are such that
$$N_j=\prod_{\nu=1}^{j} n_\nu,\quad n_j\to\infty\,(j\to\infty),
\quad \log n_{j+1}/\log N_j\to0,\,\ve_j\to0\,(j\to\infty).$$
(We consider that the empty product for $j=0$ is equal to $1$.)
Let also all $n_j$ be large enough so that for any $j$ the discs
$B(x,n_j^{-1/2-2\ve_j})$, $x\in A_j$, are mutually
disjoint and contained in $B(0,1)$; this is possible by (2.2). Denote
$$\delta_j=n_j^{-1/2-2\ve_j},\quad \Delta_j=\prod_{\nu=1}^{j} \delta_j.$$
Let $E_1=\bigcup_{x\in A_1}B(x,\delta_1)$.
We then define $E_2,E_3,\dots$ by induction.
Namely, suppose that we have constructed $E_j$ which is
a union of $N_j$ disjoint closed discs $B_i$ of radius $\Delta_j$ each.
Then $E_{j+1}$ is obtained from $E_j$ by replacing each $B_i$ by
the image of $\bigcup_{x\in A_{j+1}}B(x,\delta_{j+1})$ under the unique affine mapping
which takes $B(0,1)$ to $B_i$ and preserves direction of vectors.
We then let $E= \bigcap_{j=1}^\infty E_j$.
The verification of properties $\dim(E)=2$ and $|\Delta_X(E)|=0$ is exactly
as in the proof of Theorem 1.
\enddemo

\bigskip

{\bf Acknowledgements.} Part of this work was completed while the first author
was a PIMS Distinguished Chair at the University of British Columbia,
We also acknowledge the support of NSERC under grant 22R80520.

\bigskip

\bigskip
\centerline{REFERENCES}
\bigskip

[B94] J.~Bourgain, Hausdorff dimension and distance sets,
Israel J. Math. {\bf 87} (1994), 193--201.

[BBKM02] V.V.~Beresnevich, V.I.~Bernik, D.Y.~Kleinbock, and G.A.~Margulis,
Metric Diophantine approximation: the Khintchine---Groshev theorem
for nondegenerate manifolds, Moscow Mathematical Journal {\bf 2} (2002), 203--225.

[Be02] V.~Beresnevich, A Groshev type theorem for convergence on manifolds,
Acta Math.~Hung. {\bf 94} (2002), 99--130.

[BKM01] V.~Bernik, D.~Kleinbock, and G.~Margulis, Khintchine--type theorems
on manifolds: the convergence case for standard and multiplicative version,
Int.~Math.~Research Notices No.~9 (2001), 453--486.

[Er03] M. B. Erdogan, Falconer's distance set conjecture,
preprint, 2004.


[F85] K.J. Falconer, The geometry of fractal sets, Cambridge University
Press (1985).

[F86] K.J. Falconer, On the Hausdorff dimension of distance sets,
Mathematika {\bf 32} (1986), 206--212.

[I01] A. Iosevich, Curvature, combinatorics and the Fourier transform, Notices
Amer. Math. Soc. {\bf 46} (2001), 577--583.

[I{\L}03] A.~Iosevich and I.~{\L}aba, Distance sets of well-distributed planar
point sets, Discrete Comput.~Geometry {\bf 31} (2004), 243--250.

[I{\L}04] A.~Iosevich and I.~{\L}aba, $K$-distance sets, Falconer conjecture and
discrete analogs, preprint, 2003.

[K{\L}04] S.~Konyagin and I. {\L}aba, Distance sets of well distributed planar
sets for polygonal norms, preprint, 2004.


[KM98] D.Y.~Kleinbock and G.A.~Margulis, Flows on homogeneous spaces and
Diophantine approximation on manifolds, Ann.~Math. (2) {\bf 148} (1998),
339--360.

[M87] P. Mattila, Spherical averages of Fourier transforms of measures with
finite energy; dimension of intersections and distance sets, Mathematica
{\bf 34} (1987), 207--228.


[W99] T.~Wolff, Decay of circular means of Fourier transforms of measures,
Int. Math. Res. Notices {\bf 10} (1999), 547--567.

\bigskip

\enddocument